\documentclass{amsart}
\usepackage{amsfonts,amsthm,euscript,epsfig,color}

\newtheoremstyle{my}{1.5em}{0.5em}{\em}{}{\sc}{:}{0.5em}{}

\theoremstyle{my}

\newtheorem{thm}{Theorem}
\newtheorem{theorem}[thm]{Theorem}
\newtheorem*{theorem*}{Theorem}

\newtheorem*{corollary*}{Corollary}

\newtheorem{conjecture}[thm]{Conjecture}
\newtheorem*{conjecture*}{Conjecture}

\newtheorem*{question*}{Question}

\newtheorem*{definitions*}{Definitions}

\newtheorem*{rem*}{Remark}

\newtheorem*{remark*}{Remark}

\newtheorem*{remarks*}{Remarks}
\newtheorem*{example*}{Example}
\newtheorem{example}[thm]{Example}
\newtheorem*{examples*}{Examples}

\newcommand{\Z}{\mathbb{Z}}
\newcommand{\Q}{\mathbb{Q}}
\newcommand{\C}{\mathbb{C}}

\newcommand{\iso}{\cong}           
\newcommand{\htp}{\simeq}          

\newcommand{\CP}[1]{\C {\mathrm P}^{#1}}

\renewcommand{\o}{\omega}

\newcommand{\A}{\EuScript{A}}
\newcommand{\B}{\EuScript{B}}
\newcommand{\CC}{\EuScript{C}}
\newcommand{\D}{\EuScript{D}}
\newcommand{\MM}{\EuScript{M}}

\newcommand{\QQ}{\EuScript{Q}}
\renewcommand{\SS}{\EuScript{S}}
\newcommand{\F}{\EuScript{F}}
\newcommand{\K}{\mathbb{K}}
\newcommand{\ZZ}{\mathcal{Z}}

\title{Symplectic homology as Hochschild homology}
\author{Paul Seidel}
\begin{document}
\maketitle

\section{Introduction}

In the wake of Donaldson's pioneering work \cite{donaldson02}, Picard-Lefschetz
theory has been extended from its original context in algebraic geometry to (a
very large class of) symplectic manifolds. Informally speaking, one can view
the theory as analogous to Kirby calculus: one of its basic insights is that
one can give a (non-unique) presentation of a symplectic manifold, in terms of
a symplectic hypersurface and a collection of Lagrangian spheres (vanishing
cycles) in it. This is particularly impressive in the four-dimensional case,
since the resulting data are easy to encode combinatorially; but the formalism
works just as well in higher dimensions. These kinds of presentations are
instructive and useful in some respects, but hard to work with in others. For
instance, it is not obvious how to recover the known symplectic invariants,
such as Gromov-Witten invariants, from vanishing cycle data. In these notes, we
ask a simpler version of this question, regarding one of the basic invariants
of symplectic manifolds with boundary, namely symplectic homology as defined by
Viterbo \cite{viterbo97a} (a closely related construction is due to
Cieliebak-Floer-Hofer \cite{cieliebak-floer-hofer95}).

There are several good surveys on symplectic homology and its basic properties,
for instance \cite{oancea04b,weber06}. Very briefly, it is an invariant
(deformation invariant) of Liouville domains. Here, by a Liouville domain
we mean a compact manifold with boundary $E = E^{2d}$, equipped with a one-form
$\theta$ such that $\o = d\theta$ is symplectic, and the dual Liouville vector
field $Z$ (defined by $i_Z\o = \theta$) points strictly outwards along the
boundary. We will also assume that $c_1(E) = 0$, and in fact we want to choose
a preferred trivialization of the canonical bundle $K_E =
\lambda^{top}_{\C}(TE)$, which turns $E$ into the symplectic counterpart of an
affine Calabi-Yau variety (this is not really necessary, but it makes some
aspects more intuitive). Finally, we fix a coefficient field $\K$, which will
be used in all Floer homology type constructions. In this setup, symplectic
homology $SH_*(E)$ is a $\Z$-graded $\K$-vector space (not necessarily
finite-dimensional, not even in a fixed degree). It comes with a natural
homomorphism
\begin{equation} \label{eq:weinstein-map}
SH_*(E) \longrightarrow H_{*+d}(E;\K),
\end{equation}
which is important for applications to the Weinstein conjecture on Reeb orbits
(failure of this to be an isomorphism indicates existence of at least one
periodic Reeb orbit for every possible choice of contact one-form on $\partial
E$). The prototypical example is where $E = DT^*\!N$ is the ball cotangent
bundle of a closed oriented $d$-manifold $N$. In that case,
\begin{equation} \label{eq:loop}
SH_*(DT^*\!N) \iso H^{-*}({\EuScript L}N;\K),
\end{equation}
where ${\EuScript L}N$ is the free loop space. With respect to this
isomorphism, \eqref{eq:weinstein-map} is restriction to constant loops combined
with Poincar{\'e} duality. \eqref{eq:loop} implicitly fixes all the conventions
used in the present paper (homology versus cohomology, the grading, and the
inclusion of non-contractible loops).

Let $\pi: E \rightarrow D$ be an (exact) Lefschetz fibration over a closed disc $D$.
For the moment, the only relevant properties of such fibrations are that each
regular fibre is a Liouville domain, and that the total space becomes
a manifold of the same kind after some minor manipulations. As usual, we fix a
trivialization of $K_E$, and (in a slight abuse of notation) denote by
$SH_*(E)$ the resulting symplectic homology. Choose a base point $\ast \in
\partial D$, and let $M = E_\ast$ be the fibre over that point. Now fix a
distinguished basis of vanishing paths leading from $\ast$ to the critical
points of $\pi$, and let $(L_1,\dots,L_m)$ be the resulting basis of vanishing
cycles, which is an ordered collection of Lagrangian spheres in $M$. Each $L_j$
gives rise to an object of the Fukaya category $\F(M)$, which (to emphasize its
more algebraic role) we denote by $X_j$. Let $\B \subset \F(M)$ be the full
$A_\infty$-subcategory with objects $(X_1,\dots,X_m)$, and $\A \subset \B$ its
directed subcategory. Extend the morphism spaces in $\B$ by introducing a
formal variable $t$ of degree $2$, which yields another $A_\infty$-category
$\B[[t]]$. Next, consider the $A_\infty$-subcategory $\CC = \A \oplus t\B[[t]]
\subset \B[[t]]$, in which the constant ($t^0$) term is constrained to lie in
$\A$. Finally, turn $\CC$ into an obstructed (or curved) $A_\infty$-category,
by switching on a $\mu^0$ term which is $t$ times the identity; and denote this
gadget by $\D$.

\begin{conjecture} \label{th:main}
The Hochschild homology of $\D$ is the symplectic homology of the total space
$E$:
\begin{equation} \label{eq:main}
SH_*(E) \iso HH_*(\D).
\end{equation}
\end{conjecture}

This conjecture is the main point of these notes, and we will try to illuminate
it from various perspectives. For the moment, a few simple checks may suffice.
In the trivial case when there are no vanishing cycles at all, both $\D$ and
its Hochschild homology vanish; but on the other hand, the total space is $E =
D \times M$, whose symplectic homology is zero (by the K{\"u}nneth formula
\cite{oancea04}, for instance). More generally, suppose that we attach a
Weinstein $(d-1)$-handle to the boundary of $M$, leaving the vanishing cycles,
hence the category $\D$, unchanged. On the level of the total space, this
results in a subcritical handle attachment, which does not affect symplectic
homology \cite{cieliebak02}. These examples are meant to address an obvious
concern, namely the fact that the ordinary homology of $E$ does not enter into
(indeed, cannot be reconstructed from) $\D$. Still, going beyond such
degenerate cases, it is by no means clear why $HH_*(\D)$ should be an invariant
of $E$, and independent of the particular Lefschetz fibration. In fact, it may
be interesting to look for a direct proof of this, not based on the expected
relation with symplectic homology (or alternatively, if one wants to be
pessimistic, this might be a good way to find a counterexample).

The plan for the rest of these notes will be as follows. The first few sections
cover preliminaries, both geometric and algebraic. We then define more
carefully the algebraic objects involved in Conjecture \ref{th:main}, and
mention some example computations. Following that, we give a speculative
geometric interpretation of the basic spectral sequence with target $HH_*(\D)$.
Finally, going beyond Hochschild homology, we take a look at the main new
object $\D$ itself, from a mostly algebraic perspective. The material in these
last sections is increasingly tentative: {\em reader, beware!} To keep the
discussion focused, we have excluded a number of related topics. For instance,
it is natural to compare $S^1$-equivariant symplectic homology to cyclic
homology (in their various respective versions), but we will not explicitly
address this, even though cyclic homology appears briefly in Section
\ref{sec:periodic}. Along related lines, symplectic homology carries a rich
structure of homology operations, and one can ask for their algebraic
counterparts in Hochschild homology (for anyone interested in this,
\cite{costello05} is a good place to start). Finally, there are analogues of
symplectic homology in the context of Lagrangian Floer homology (for Lagrangian
submanifolds with boundary); again, these may be glimpsed in Section
\ref{sec:iterate}, but will not appear directly.

{\em Acknowledgements.} Several years ago, Donaldson and Eliashberg
independently suggested that I should look at constructions somewhat similar to
$HH_*(\D)$ (I hereby apologize for being late in responding to their ideas!).
In Donaldson's case, this was motivated by thinking about the boundary
monodromy of the Lefschetz fibration and its iterates (a point of view which
will be adopted in Section \ref{sec:periodic}). Eliashberg was interested in
the change of contact homology under Legendrian surgery (this is currently
being pursued in joint work of Bourgeois, Ekholm and Eliashberg); the relation
between the two topics is established by ongoing work of Oancea and Bourgeois,
who constructed a long exact sequence relating symplectic and contact homology.
Numerous discussions of symplectic homology with Ivan Smith, and of cyclic
homology with Kevin Costello, have been enormously helpful. I'd also like to
thank the organizers of the 2005 AMS Summer Institute in Algebraic Geometry for
allowing me to present some rather half-baked ideas. This research was
partially funded by NSF grant DMS-0405516.

\section{Lefschetz fibrations and categories}

First, we need to clarify the notion of Lefschetz fibration involved. This is
easiest to explain in the case when there are no critical points. Then, what we
want to have is a differentiable fibre bundle $\pi: E \rightarrow D$, whose
fibres are compact manifolds with boundary (which of course means that $E$ has
codimension $2$ corners), together with a one-form $\theta$ whose restriction
to each fibre gives it the structure of a Liouville domain. According to the standard
theory of symplectic fibrations, $\o = d\theta$ determines a symplectic (in
fact, Hamiltonian) connection. This will not in general have well-defined
parallel transport, because the integral flow lines can hit the boundary of the
fibres. However, this deficiency can be easily corrected by deforming $\theta$
in an appropriate way; see \cite[Section 6]{khovanov-seidel98}. 
Assume from now on that this has been done. Then, by adding a large
multiple of a suitable one-form from the base, we can achieve that $\o$ itself
becomes symplectic, and that the Liouville field dual to $\theta$ points
strictly outwards along all boundary faces; compare \cite[Section 15b]{seidel04}. After that, rounding off the
corners turns $E$ itself into a Liouville domain. The outcome is independent
of all the details up to deformation. The Lefschetz case is largely the same,
except that we allow $\pi$ to have finitely many critical points, lying in the
interior of $E$, and locally (symplectically) modelled on nondegenerate
singular points of holomorphic functions on a K{\"a}hler manifold. For
simplicity, we require that there should be at most one such point in each
fibre. The previous discussion carries over with some small modifications; in
particular, the outcome justifies the notation $SH_*(E)$ used in Conjecture
\ref{th:main}.

As before, we fix some $\ast \in
\partial D$, which by definition is a regular value of $\pi$, and set $M =
E_\ast$. Recall that in our context, $E$ always comes with a trivialization of
its canonical bundle, which is then inherited by the fibre because $K_M \iso
K_E|M$. Given that, one can define the Fukaya category $\F(M)$, which is an
$A_\infty$-category linear over $\K$. Objects are (exact) closed Lagrangian
submanifolds decorated with some minor additional data (grading, $Spin$
structure). Roughly speaking, the morphism spaces and their differentials
$\mu^1$ are given by the cochain complexes underlying Lagrangian Floer
cohomology theory, $hom_{\F(M)}(L_0,L_1) = CF^*(L_0,L_1)$; and the higher order
compositions $\mu^d$ count pseudo-holomorphic polygons. There are various ways
to implement the details, which differ in how they deal with technical issues
such as non-transversally intersecting Lagrangian submanifolds, but the
resulting $A_\infty$-categories are all quasi-isomorphic to each other. We will
assume that things have been arranged in such a way that $\F(M)$ is strictly
unital (has cochain level identity morphisms). This may not be the case with
the most common definitions, which only produce cohomology level identities,
but that can always be amended by passing to a quasi-isomorphic
$A_\infty$-structure (without changing the morphism spaces themselves; this is
a general algebraic fact, see \cite[Chapter 3]{lefevre} or \cite[Section 2]{seidel04}). Suppose from now on that our Lefschetz fibration has
well-defined parallel transport maps. Choose a distinguished basis of paths
$(\gamma_1,\dots,\gamma_m)$ going from $\ast$ to the critical values of $\pi$.
Each $\gamma_j$ gives rise to a Lefschetz thimble $\Delta_j$, which is a
Lagrangian ball in $E$ fibered over $\gamma_j$. Its boundary, the vanishing
cycle $L_j = \partial\Delta_j$, is a Lagrangian sphere in $M$
\cite[Section 1.3]{exact}. We equip the
latter with the necessary additional data to make them into objects $X_j$ of
$\F(M)$ (the grading may be chosen arbitrarily; the $Spin$ structure is
inherited from $\Delta_j$, hence will be the nontrivial one for $d = 2$).
Denote by $\B \subset \F(M)$ the associated full $A_\infty$-subcategory. To
remember the ordering of the objects $X_j$, we also consider the directed
subcategory $\A \subset \B$, in which the morphism spaces are by definition
\begin{equation}
hom_\A(X_j,X_k) = \begin{cases} hom_\B(X_j,X_k) & j<k, \\ \K e_j & j = k,
\\ 0 & j>k; \end{cases}
\end{equation}
$e_j \in hom_\B(X_j,X_j)$ being the strict identity. $\A$ is probably the more
accessible of the two structures. For instance, if $d = 2$, it can be
determined in a simple combinatorial way, by counting immersed polygons on $M$
with boundary in the $L_j$ \cite[Section 13]{seidel04}. Beyond that, there are other methods for computing
both $\A$ and $\B$, which work in higher dimensions too, but they are
considerably less elementary.

Given an $A_\infty$-category $\ZZ$, one can define its derived category
$D(\ZZ)$ as follows \cite{keller}. Let $mod(\ZZ)$ be the $A_\infty$-category (in fact
differential graded category, since all compositions of order $>2$ vanish) of
right $\ZZ$-modules. For our present purpose, it is convenient to assume that
$\ZZ$ is strictly unital, and to use strictly unital modules. There is a
natural $A_\infty$-functor $\ZZ \rightarrow mod(\ZZ)$, the Yoneda embedding,
which is full and faithful on the cohomology level; see for instance \cite[Section 2]{seidel04}. 
Moreover, $mod(\ZZ)$ is a
triangulated $A_\infty$-category, which means that it is closed under shifts
and forming mapping cones (unlike the corresponding notion in classical
homological algebra, being triangulated is a property rather than an additional
structure). Then, $D(\ZZ)$ is the smallest triangulated subcategory of
$mod(\ZZ)$ which contains the image of the Yoneda embedding. Note that in the
case of a directed $A_\infty$-category, such as $\A$, the derived category will
be (quasi-equivalent to) the category of all finite-dimensional
$A_\infty$-modules. Our interest in derived categories comes from the following
result:

\begin{theorem} \label{th:derived}
Up to quasi-equivalence, $D(\A)$ and $D(\B)$ are independent of the choice of
vanishing cycles, hence invariants of the Lefschetz fibration.
\end{theorem}

The first step in proving this is to note that a smooth isotopy of the paths
$\gamma_i$ gives rise to a Lagrangian (exact Lagrangian, to be precise) isotopy
of the vanishing cycles, which leaves the categories $\A$ and $\B$ unchanged up
to quasi-isomorphism. This is a version of the standard isotopy invariance
property of Floer theory. Next, considering only isotopy classes, there is a
simply-transitive action of the braid group $Br_m$ on the set of disinguished
bases of paths. On the level of vanishing cycles, this yields the so-called
Hurwitz moves which relate any two bases of such cycles. For instance, the
standard $i$-th generator $s_i$ of $Br_m$ gives rise to an elementary Hurwitz
move
\begin{equation} \label{eq:hurwitz}
(L_1,\dots,L_m) \longmapsto
(L_1,\dots,L_{i-1},\tau_{L_i}(L_{i+1}),L_i,L_{i+2},\dots,L_m),
\end{equation}
where $\tau_{L_i}$ is the symplectic Dehn twist along $L_i$. Note that $D(\B)$
is obviously equivalent to the triangulated subcategory of $D(\F(M))$ generated
by the $L_i$, which we denote by ${\mathcal T}$. The effect of $\tau_{L_i}$ on
objects of $\F(M)$ is well-known: it corresponds to the so-called twist functor
associated to the spherical object $L_i$ \cite[Corollary 17.17]{seidel04}. From this interpretation, it follows
that the modified basis of vanishing cycles obtained by a Hurwitz move
\eqref{eq:hurwitz} again lies in ${\mathcal T}$, and generates that category.
In other words, ${\mathcal T} \subset D(\F(M))$ is independent of the choice of
basis, which proves the part of Theorem \ref{th:derived} concerning $D(\B)$.
Note that the result is in fact slightly stronger than stated: there is a
unique object, namely ${\mathcal T}$, which is canonically equivalent to
$D(\B)$; therefore, if $\B,\B'$ are the categories arising from two choices of
bases, the functor $D(\B) \iso D(\B')$ is essentially canonical, and satisfies
the obvious composition property.

The story for $\A$ is a little more involved. On the algebraic side, we have
the theory of mutations, which gives rise to an action of $Br_m$ on
(quasi-isomorphism classes of) directed $A_\infty$-categories of length $m$.
Moreover, if $\A,\A'$ are two categories lying in the same orbit, then $D(\A)
\iso D(\A')$ (for a survey of mutations and their applications, see
\cite{gorodentsev-kuleshov04}; the generalization to $A_\infty$-categories, and
the idea of applying this to Lefschetz fibrations, are due to Kontsevich). One
can use the previously mentioned results on Dehn twists to show that if we take
$\A$, and apply the mutation corresponding to a generator of $Br_m$, the
outcome corresponds to the effect of a move \eqref{eq:hurwitz} on the vanishing
cycles. This proves the part of Theorem \ref{th:derived} concerning $D(\A)$ as
stated, which is \cite[Theorem 17.20]{seidel04}, but does not by itself establish canonical equivalences. 
To get around
this problem, one can introduce another object $\F(\pi)$, the Fukaya category
of the Lefschetz fibration $\pi$, and then prove that for each choice of basis,
there is an embedding $\A \rightarrow \F(\pi)$ which induces an equivalence of
derived categories (in other words, every basis forms a full exceptional
collection in $D(\F(\pi))$, and the directed category associated to that
collection is $\A$). This is carried out in \cite[Section 18]{seidel04}, under the
technical assumption that $char(\K) \neq 2$. The details are not terribly
relevant for our present discussion, but we should give some idea of the
geometry behind $\F(\pi)$. Objects are Lefschetz thimbles $\Delta \subset E$
with the usual additional data. When forming their morphism group, which is
again written as $CF^*(\Delta_0,\Delta_1)$, one perturbs the path $\gamma_0$
underlying $\Delta_0$ slightly by moving its endpoint in positive direction
along the boundary. This gets rid of boundary intersection points, allowing one
to apply the standard Floer-theoretic formalism. The reader should now see why
a basis $(\Delta_1,\dots,\Delta_m)$ gives rise to an exceptional collection in
$\F(\pi)$: for $i>j$, we have $CF^*(\Delta_i,\Delta_j) = 0$ because the
perturbation of $\gamma_i$ makes the two Lagrangian submanifolds disjoint
(Figure \ref{fig:disjoint}).

\begin{figure}
\begin{centering}
\begin{picture}(0,0)%
\includegraphics{perturb.pstex}%
\end{picture}%
\setlength{\unitlength}{3947sp}%
\begingroup\makeatletter\ifx\SetFigFont\undefined%
\gdef\SetFigFont#1#2#3#4#5{%
  \reset@font\fontsize{#1}{#2pt}%
  \fontfamily{#3}\fontseries{#4}\fontshape{#5}%
  \selectfont}%
\fi\endgroup%
\begin{picture}(2416,2419)(2393,-1423)
\put(2776,-736){\makebox(0,0)[lb]{\smash{{\SetFigFont{10}{12.0}{\rmdefault}{\mddefault}{\updefault}{$\gamma_1$}%
}}}}

\put(3526,-286){\parbox{7em}{$\gamma_2$ and its per\-turbed version}}
%
%
\end{picture}%
\caption{\label{fig:disjoint}}
\end{centering}
\end{figure}

\section{Global monodromy}
The invariants $D(\A)$ and $D(\B)$ are of a fairly abstract kind. For instance,
it is not clear how to use them to efficiently distinguish between Lefschetz
fibrations. It is therefore natural to ask how they relate to more
straightforward geometric invariants, such as the monodromy group (viewed as a
subgroup of the symplectic mapping class group of $M$). We will not try to
address this systematically, but we will give one (conjectural) example of such
a relation, which was suggested by Donaldson. Assume that our Lefschetz
fibration has well-defined parallel transport maps; then, by going in positive
sense around the boundary, we get an automorphism $\mu$ of $M$, called the
global monodromy. From the trivialization of $K_E$, this inherits a little bit
of additional data (a grading), which means that there is a well-defined
associated graded vector space $HF_*(\mu)$, the (fixed point) Floer homology of
$\mu$. Such Floer homology groups can be thought of as part of a TQFT-type
structure arising from Lefschetz fibrations. In particular, in our case we have
a natural map
\begin{equation} \label{eq:sections}
H_*(M;\K) = HF_*(id) \longrightarrow HF_*(\mu),
\end{equation}
defined roughly speaking by counting pseudo-holomorphic sections of $\pi: E
\rightarrow D$. By taking mapping cones on the chain level, define a relative
group $HF_*(\mu,id)$ which fits into a long exact sequence
\begin{equation}
\cdots \rightarrow H_*(M;\K) \longrightarrow HF_*(\mu) \longrightarrow
HF_*(\mu,id) \rightarrow \cdots
\end{equation}
It should be emphasized that the map \eqref{eq:sections}, and hence also the
relative group $HF_*(\mu,id)$, depend on the Lefschetz fibration, and not just
on $\mu$ and $M$.

Given any basis $(L_1,\dots,L_m)$ of vanishing cycles, one can write the global
monodromy (up to Hamiltonian isotopy) as
\begin{equation}
\mu \htp \tau_{L_1} \cdots \tau_{L_m}.
\end{equation}
By the long exact sequence from \cite{seidel00b}, this means that in principle,
the difference between the two groups in \eqref{eq:sections} can be expressed
in terms of Lagrangian Floer cohomology groups of the vanishing cycles. To make
this idea more precise, we build a chain complex
\begin{equation} \label{eq:donaldson}
\bigoplus_{\substack{n \geq 0 \\
j_0 < \cdots < j_n}} \!\!\! \big(CF^*(X_{j_n},X_{j_0}) \otimes
CF^*(X_{j_{n-1}},X_{j_n}) \otimes \cdots \otimes CF^*(X_{j_0},X_{j_1})
\big)[d+n],
\end{equation}
Graphically, the generators of this complex can be thought of as closed chains
of morphisms; more precisely, if one thinks of the critical values as arranged
on a circle, the chain of morphisms goes once around that circle, hence can be
drawn as an inscribed polygon. The differential consists of using the
$A_\infty$-products in $\F(M)$ to shorten the polygons in all possible ways,
see Figure \ref{fig:polygon}. Alternatively, to give a more concise algebraic
formulation, recall that there is a bijective correspondence between
$A_\infty$-categories with $m$ numbered objects, and $A_\infty$-algebras linear
over the semisimple ring $R = \K^m = \K e_1 \oplus \cdots \oplus \K e_m$. This
correspondence is given by taking an $A_\infty$-category $\ZZ$ and turning it
into the algebra $\bigoplus_{j,k} hom_\ZZ(X_j,X_k)$, with the left and right
action of the $e_j$ given by projection onto the various summands. With that in
mind, take $\A_+$ to be the augmentation ideal of $\A$, which is the kernel of
the obvious map $\A \rightarrow R$; and let $T(\A_+[1]) = R \oplus \A_+[1]
\oplus \A_+[1] \otimes_R \A_+[1] \oplus \cdots$ be the tensor algebra over $R$
with generators $\A_+[1]$ (this tensor algebra is actually finite-dimensional,
since the tensor product of $m$ copies of $\A_+$ vanishes by directedness). One
can then write \eqref{eq:donaldson} as
\begin{equation} \label{eq:donaldson2}
(\B[d] \otimes T(\A_+[1]))^{diag}
\end{equation}
where the tensor product is again taken over $R$, and the superscript $diag$
means that we retain only the diagonal piece $\bigoplus_i e_iQe_i$ of an
$R$-bimodule $Q$. The differential on \eqref{eq:donaldson2} is a Hochschild
type expression (we will not write it down here, but the connection will be
made explicitly in Section \ref{sec:periodic}). Donaldson's conjecture, in a
form closely related to the one in \cite{seidel00b}, is:

\begin{conjecture} \label{th:donaldson}
The cohomology of \eqref{eq:donaldson} (after grading-reversal) is isomorphic
to $HF_*(\mu,id)$.
\end{conjecture}

While this remains unproved at the moment, it seems well within reach of the
existing technology (particularly in view of ongoing work of
Wehrheim-Woodward).
\begin{figure}
\begin{centering}
\begin{picture}(0,0)%
\includegraphics{polygon.pstex}%
\end{picture}%
\setlength{\unitlength}{3947sp}%
\begingroup\makeatletter\ifx\SetFigFont\undefined%
\gdef\SetFigFont#1#2#3#4#5{%
  \reset@font\fontsize{#1}{#2pt}%
  \fontfamily{#3}\fontseries{#4}\fontshape{#5}%
  \selectfont}%
\fi\endgroup%
\begin{picture}(2673,2505)(2101,-1469)
\put(2896, 14){\begin{minipage}{6em}$\mu^2(a_1,a_4) \in$\newline ${\,\,}
CF^*(L_5,L_2)$\end{minipage}}
\put(1356,-661){\makebox(0,0)[lb]{\smash{{\SetFigFont{10}{12}{\rmdefault}{\mddefault}{\updefault}{\color[rgb]{0,0,0}$a_1 \in CF^*(L_1,L_2)$}%
}}}}
\put(1651,614){\makebox(0,0)[lb]{\smash{{\SetFigFont{10}{12}{\rmdefault}{\mddefault}{\updefault}{\color[rgb]{0,0,0}$a_2 \in CF^*(L_2,L_3)$}%
}}}}
\put(4131,-361){\makebox(0,0)[lb]{\smash{{\SetFigFont{10}{12}{\rmdefault}{\mddefault}{\updefault}{\color[rgb]{0,0,0}$a_3 \in CF^*(L_3,L_5)$}%
}}}}
\put(3031,-1411){\makebox(0,0)[lb]{\smash{{\SetFigFont{10}{12}{\rmdefault}{\mddefault}{\updefault}{\color[rgb]{0,0,0}$a_4 \in CF^*(L_5,L_1)$}%
}}}}
\end{picture}%
\caption{\label{fig:polygon}}
\end{centering}
\end{figure}

\section{The Serre functor}
There is one piece of the puzzle which we haven't mentioned so far. Suppose
that our Lefschetz fibration is such that parallel transport is well-defined
and the total space is symplectic. In addition, we want the fibration to be
symplectically locally trivial near $\partial D$, which can also be achieved by
a suitable deformation of $\theta$. Take a circle $\lambda$ in $int(D)$ which
runs parallel to the boundary, and consider the negative Dehn twist (in the
classical sense) $t_{\lambda}^{-1}$ along it. Under the assumptions we have
made, this lifts to a symplectic automorphism of $E$, which we denote by
$\sigma$. By definition, $\sigma|M = \mu$, while for any $z$ which is
sufficiently far from $\partial D$ we have $\sigma|E_z = id$. For general
reasons, $\sigma$ gives rise to an auto-equivalence of the category $\F(\pi)$
as well as its derived category, hence also of $D(\A)$. We denote all of these
by $\sigma_*$. It is maybe instructive to reformulate the construction in a
slightly different way. As mentioned before, the braid group $Br_m$ acts on the
set of all possible choices of vanishing paths, in a way which preserves the
derived category of $\A$. More precisely, if one applies an element of the
braid group to a given set of vanishing paths, the result is another directed
Fukaya category $\A'$ such that
\begin{equation} \label{eq:deq}
D(\A) \iso D(\A');
\end{equation}
to be more precise, we want to use the preferred equivalence which comes from
comparing both categories with $D(\F(\pi))$. Now consider the case when our
element is the distinguished positive generator of the center, $b \in Z(Br_m)
\iso \Z$. The action of this on vanishing paths is just given by
$t_{\lambda}^{-1}$; on the level of vanishing cycles, this means that one
applies $\mu$ to all cycles simultaneously. This of course does not change the
associated Fukaya category, so that actually $D(\A) = D(\A')$. In view of this,
the equivalence \eqref{eq:deq} turns into an automorphism of $D(\A)$, which is
$\sigma_*$. What is interesting for us is that due to the geometry of $\sigma$,
there is a canonical natural transformation
\begin{equation} \label{eq:geom-tr}
\sigma_* \longrightarrow id.
\end{equation}
This is best visible in terms of $\F(\pi)$, where one can indeed arrange that
for any two Lefschetz thimbles $(\Delta_0,\Delta_1)$, $CF^*(\Delta_0,\Delta_1)$
is a subcomplex of $CF^*(\sigma(\Delta_0),\Delta_1)$. This in particular yields
a canonical element in $HF^*(\sigma(\Delta),\Delta)$, namely the image of the
identity in $HF^*(\Delta,\Delta)$, which is the cohomology level natural
transformation induced by \eqref{eq:geom-tr}.

It is an old idea (due to Kontsevich, I believe) that $\sigma_*$ should in fact
be the Serre functor $\SS$ of $D(\A)$, up to a shift by $d$. Recall that in the
classical categorical context, Serre functors $S$ are defined by the existence
of a natural isomorphism $Hom(Y_0,S(Y_1)) \iso Hom(Y_1,Y_0)^\vee$. There is a
refinement of this on the $A_\infty$-level, which determines $\SS$ up to
isomorphism. The details don't really matter here, but we do want to make the
point that $\SS$ is intrinsic to $D(\A)$, and does not depend on any additional
choices. While Kontsevich's conjecture has not yet been completely proved, one
can show that $\sigma_* \iso \SS$ as far as the action on objects is concerned
(this is a consequence of the basic relation between Hurwitz moves on vanishing
cycles and mutations, together with a purely algebraic argument relating
mutation by $b$ to the Serre functor \cite[Assertion 4.2]{bondal89}). In view
of this, we feel free to assume that the functor isomorphism is indeed true,
the outcome being that $D(\A)$ comes with a preferred natural transformation
$\SS \rightarrow id$ of degree $d$.

Remarkably, there is another construction of such a natural transformation,
which is much more algebraic, and a priori has nothing to do with the geometry
of Lefschetz fibrations. Recall first that, given a finite-dimensional
$A_\infty$-bimodule $\QQ$ over $\A$, one can define a differential graded
functor from $D(\A)$ (thought of as category of finite-dimensional
$A_\infty$-modules) to itself, whose effect on objects is
\begin{equation} \label{eq:derived-tensor}
\MM \longmapsto \MM \otimes_\A \QQ = \MM \otimes T(\A_+[1]) \otimes \QQ.
\end{equation}
This is just the definition of the tensor product in the $A_\infty$-context; on
the right hand side, the tensor products are over $R$, and the resulting graded
vector space comes with a canonical $A_\infty$-module structure. (To motivate
the definition, look at the classical context, where we have an algebra $A$, a
bimodule $Q$, and a module $M$; there, the tensor product in
\eqref{eq:derived-tensor} would be called a derived one, since it's obtained by
taking $M \otimes_A Q = M \otimes_A A \otimes_A Q$ and then replacing the $A$
in the middle with its reduced bar resolution.) For instance, if one takes $\QQ
= \A$ to be the diagonal bimodule, the resulting functor is isomorphic to the
identity. More interestingly, one can take its dual $\QQ = \A^\vee$, and then
\eqref{eq:derived-tensor} is isomorphic to $\SS$. In our situation, the
embedding $\A \subset \B$ gives rise to a short exact sequence of bimodules
\begin{equation} \label{eq:short}
0 \rightarrow \A \longrightarrow \B \longrightarrow \QQ = \B/\A \rightarrow 0.
\end{equation}
By inverting the quasi-isomorphism $\QQ \iso Cone(\A \rightarrow \B)$ one gets
a boundary map $\QQ \rightarrow \A[1]$ (in the derived category of
$A_\infty$-bimodules, this completes \eqref{eq:short} to an exact triangle).
Now $\B$, being part of the Fukaya category of $M$, has a weak version of
cyclic symmetry \cite{tradler01}, and this implies that $\QQ \iso \A^\vee[1-d]$
canonically. With that in mind, our boundary map becomes a bimodule
homomorphism
\begin{equation} \label{eq:alg-tr}
\A^\vee \longrightarrow \A[d];
\end{equation}
on the level of functors, this induces a natural transformation $\SS
\rightarrow id$ of degree $d$. It is natural to propose the following

\begin{conjecture} \label{th:agree}
The natural transformations obtained from \eqref{eq:geom-tr} and
\eqref{eq:alg-tr} agree.
\end{conjecture}

\begin{example} \label{th:mirror}
An interesting class of examples is provided by mirror symmetry. We briefly
recall the setup, using classical derived categories for simplicity (and of
course, setting $\K = \C$). Let $X$ be a smooth $d$-dimensional Fano toric
variety, and $Y \subset X$ the toric anticanonical divisor, which is simply the
union of all codimension $1$ torus orbits. Consider the derived categories of
coherent sheaves $D^b(Y)$ and $D^b(X)$, where in the first case we limit
ourselves to perfect complexes, to avoid problems arising from singularities.
$D^b(Y)$ is a Calabi-Yau category, meaning that its Serre functor is a shift
$[d-1]$. In the case of $D^b(X)$, the Serre functor is $S = K_X[d] \otimes -$,
which means that the section $s$ defining $Y$ gives rise to a natural
transformation $S[-d] \rightarrow id$. In fact, for any object $F \in D^b(X)$
we have an exact triangle
\begin{equation} \label{eq:resolution}
 F \longrightarrow i_*i^*F \longrightarrow S(F)[1-d] \longrightarrow F[1],
\end{equation}
where $i: Y \rightarrow X$ is the inclusion, and the last map in that triangle
is precisely the natural transformation we just discussed. The mirror of $X$ is
a Landau-Ginzburg theory with superpotential $W: (\C^*)^d \rightarrow \C$. If
we cut out a compact piece of this to make a Lefschetz fibration $\pi: E
\rightarrow D$, then homological mirror symmetry predicts that $D^b(X)$ should
be $H^0(D(\A))$, and similarly $D^b(Y)$ should become isomorphic to
$H^0(D(\B))$ after Karoubi completion on both sides. It seems natural to think
that $i_*i^*$ would then correspond to the endofunctor of $H(D(\A))$ induced by
the $\A$-bimodule $\B$, so that \eqref{eq:resolution} would turn into a weak
version of \eqref{eq:short}. This means that the natural transformation from
Conjecture \ref{th:agree} would correspond to $s \in H^0(X,K_X^{-1})$.
Inspection of known cases, such as $X = \CP{3}$, seems to bear out this idea.
\end{example}

\section{Hochschild homology}
We should first explain, in a little more detail than before, the definitions
of $\CC$ and $\D$. Start with the pair of categories $\A \subset \B$. $\CC$ has
the same objects; and morphisms in it are formal power series
\begin{equation}
x = x_0 + t x_1  + \cdots \in hom_\B(X_j,X_k)[[t]]
\end{equation}
where the variable $t$ has degree $2$, subject to the condition that the
leading order term $x_0$ must lie in $hom_\A(X_j,X_k)$. We equip this with the
$A_\infty$-operations obtained by $t$-linearly extending those in $\B$. $\D$ is
the same as $\CC$, except that it carries an extra curvature term, given by
\begin{equation} \label{eq:m0}
\mu^0_\D = t e_j \in hom_\CC(X_j,X_j) = hom_\D(X_j,X_j)
\end{equation}
for all $j$. Note that in spite of the obvious $t$-linearity of the composition
maps, the ground field here is still considered to be $\K$. However, we do want
to take into account the (complete) $t$-adic topology on $\CC$ and $\D$, and
that will affect all associated constructions.

The Hochschild homology $HH_*(\D)$ is defined through the reduced version of
the classical cyclic bar complex $\bar{C}^* = \bar{C}^*(\D)$. Written in terms
of $R$-modules, this is the $t$-adic completion of
\begin{equation} \label{eq:reduced-hochschild}
\big(\D \otimes T(\D_+[1])\big)^{diag},
\end{equation}
where the augmentation ideal is $\D_+ = \A_+ \oplus t\B[[t]] \subset \D$. To
make this more explicit, one can separate out the various powers of $t$. The
outcome is that $\bar{C}^*$ is the direct product of pieces
\begin{equation} \label{eq:pieces}
\begin{aligned}
& \big(\A \otimes T(\A_+[1]) \otimes t^{i_k} \B[1] \otimes T(\A_+[1]) \otimes
\cdots \otimes t^{i_1}
\B[1] \otimes T(\A_+[1])\big)^{diag}, \\
& \big( t^{i_k} \B \otimes T(\A_+[1]) \otimes t^{i_{k-1}}\B[1] \otimes \cdots
\otimes t^{i_1}\B[1] \otimes T(\A_+[1])\big)^{diag}
\end{aligned}
\end{equation}
ranging over all $k$ and $i_1,\dots,i_k \geq 1$. To represent things
graphically in analogy with Figure \ref{fig:polygon}, put the critical values
of $\pi$ on a circle, and mark a point in its center. Then, a generator of
\eqref{eq:pieces} can be drawn as a closed chain of morphisms winding $i_1 +
\cdots + i_k$ times (clockwise) around the central point, see Figure
\ref{fig:polygon2}. Note the slight asymmetry: $e_j \in \A$ can only occur as
the last (leftmost) element of the chain (in contrast, $t^ie_j$ with $i>0$ can
occur anywhere).
\begin{figure}
\begin{centering}
\begin{picture}(0,0)%
\includegraphics{polygon2.pstex}%
\end{picture}%
\setlength{\unitlength}{3947sp}%
\begingroup\makeatletter\ifx\SetFigFont\undefined%
\gdef\SetFigFont#1#2#3#4#5{%
  \reset@font\fontsize{#1}{#2pt}%
  \fontfamily{#3}\fontseries{#4}\fontshape{#5}%
  \selectfont}%
\fi\endgroup%
\begin{picture}(3399,4090)(1801,-2819)
\put(1401,-2261){
\begin{minipage}{30em} $a_5 \otimes \cdots \otimes a_1 \in
\big(t^2\B \otimes t \B[1] \otimes \A_+[1] \otimes t \B[1] \otimes \A_+[1]\big)^{diag},$ \\
consisting of: \\[1em]
$a_1 \in CF^{*+1}(L_1,L_3), \;\; a_2 \in CF^{*-1}(L_3,L_3)$, \\
$a_3 \in CF^{*+1}(L_3,L_4), \;\; a_4 \in CF^{*-1}(L_4,L_2)$, \\ $a_5 \in
CF^{*-4}(L_2,L_1).$
\end{minipage}
}
\put(3526,1139){\makebox(0,0)[lb]{\smash{{\SetFigFont{10}{12.0}{\rmdefault}{\mddefault}{\updefault}{\color[rgb]{0,0,0}$L_3$}%
}}}}
\put(4801,314){\makebox(0,0)[lb]{\smash{{\SetFigFont{10}{12.0}{\rmdefault}{\mddefault}{\updefault}{\color[rgb]{0,0,0}$L_4$}%
}}}}
\put(2101,314){\makebox(0,0)[lb]{\smash{{\SetFigFont{10}{12.0}{\rmdefault}{\mddefault}{\updefault}{\color[rgb]{0,0,0}$L_2$}%
}}}}
\put(4351,-1336){\makebox(0,0)[lb]{\smash{{\SetFigFont{10}{12.0}{\rmdefault}{\mddefault}{\updefault}{\color[rgb]{0,0,0}$L_5$}%
}}}}
\put(2626,-1336){\makebox(0,0)[lb]{\smash{{\SetFigFont{10}{12.0}{\rmdefault}{\mddefault}{\updefault}{\color[rgb]{0,0,0}$L_1$}%
}}}}
\end{picture}%
\caption{\label{fig:polygon2}}
\end{centering}
\end{figure}

The Hochschild differential $b: \bar{C}^* \rightarrow \bar{C}^{*+1}$ has three
constituents. The most familiar one is the standard bar differential induced by
the $A_\infty$-structure of $\B$, given by adding up
\begin{equation} \label{eq:d1}
 x_n \otimes \cdots \otimes x_0 \longmapsto (-1)^{\S}
 x_n \otimes \cdots \otimes
 \mu^j_\B(x_{i+j-1},\dots,x_i) \otimes x_{i-1} \otimes \cdots \otimes x_0
\end{equation}
over all $i \geq 0$ and $1 \leq j \leq n-i+1$. Here, $\S =
\|x_0\|+\cdots+\|x_{i-1}\|$, where $\|x_k\| = |x_k|-1$ denotes the reduced
grading. Next we have the additional terms which are specific of Hochschild
type theories, involving a cylic permutation of the factors:
\begin{equation} \label{eq:d2}
 x_n \otimes \cdots \otimes x_0 \longmapsto (-1)^\ast
 \mu_\B^d(x_{i-1},\dots,x_0,x_n,\dots,x_{i+j}) \otimes x_{i+j-1} \otimes
 \cdots \otimes x_i
\end{equation}
for all $i>0$ and $0 \leq j \leq n-i$. Here, $\ast =
(\|x_0\|+\cdots+\|x_{i-1}\|)(\|x_i\|+\cdots+\|x_n\|) + (\|x_i\| + \cdots +
\|x_{i+j-1}\|)$. Graphically speaking, these can both be viewed as shortening
operations similar to those in Figure \ref{fig:polygon}. Finally, reflecting
the presence of the curvature \eqref{eq:m0}, we have another term which
consists of inserting that quantity in all possible places except the first
one. Concretely, this is the sum of
\begin{equation} \label{eq:d3}
 x_n \otimes \cdots \otimes x_0 \longmapsto
 (-1)^{\S} x_n \otimes \cdots \otimes x_i \otimes
 te_{k_i} \otimes x_{i-1} \cdots \otimes x_0
\end{equation}
over $0 \leq i \leq n$, with the same sign as before. The indices $k_i$ are
such that the object $X_{k_i}$ is the source of the morphism $x_i$ (and target
of $x_{i-1}$, or of $x_n$ if $i = 0$). Finally, note that $(\bar{C}^*,b)$ as we
have written it, is a cohomological complex; to get Hochschild homology, the
sign of the grading needs to be reversed, and we will do so tacitly whenever
this issue arises.

By definition, our complex comes with a complete decreasing filtration $F^* =
F^*\bar{C}^*$ by $t$-adic weights ($F^p$ is the product of all
\eqref{eq:pieces} with $i_1 + \cdots + i_k \geq p$). All terms in the
differential preserve the weight, except for \eqref{eq:d3} which raises it by
one. Therefore, $Gr(F^*) = \prod_p F^p/F^{p+1}$ with its induced differential
is just the reduced cyclic bar complex of $\CC$. In other words, if we consider
the associated spectral sequence, whose $E_\infty$ term is the induced
filtration of $HH_*(\D)$, then the starting page can be identified with
\begin{equation} \label{eq:e1}
E^1_{pq} = \begin{cases} HH_{p+q}(\CC)^{-p} & p \leq 0, \\ 0 &
\text{otherwise,} \end{cases}
\end{equation}
where the superscripts denote $t$-weights. The negative index $-p$ appears
because this is formulated as a homology spectral sequence; the $d^1$
differential, which is the map induced by \eqref{eq:d3}, has bidegree $(-1,0)$.
The last nontrivial column is always $E^1_{00} = R$, $E^1_{0q} = 0$ for $q \neq
0$, which means that the edge homomorphism of the spectral sequence is a map
$HH_0(\D) \rightarrow H_0(\bar{C}^*/F^1) = R$. The suggested geometric
interpretation is that this would be the map $SH_0(E) \rightarrow H_d(E;\K)
\rightarrow H_d(E,M;\K)$, where the latter group is identified with $R$ by
taking the Lefschetz thimbles as a basis.

\begin{example} \label{th:extension}
Given any $A_\infty$-algebra $\A$ and bimodule $\QQ$, one can define the
extension algebra $\B = \A \oplus \QQ$. Suppose from now on that our $\B$ is
indeed of this form (this implies that the short exact sequence
\eqref{eq:short} splits, hence that the natural transformation coming from
\eqref{eq:alg-tr} is zero). Then $\CC \iso \B'[[t]]$, where $\B' = \A \oplus
\QQ[2]$. Using a slight variation of the K{\"u}nneth formula for Hochschild
homology \cite[Theorem 4.2.5]{loday}, the starting term \eqref{eq:e1} can be
written as $HH_*(\CC) \iso HH_*(\B')[[t]] \oplus HH_*(\B')[[t]]dt$, where $dt$
is a formal symbol of (homological) degree $-1$.
The differential $d^1$ is wedge product with $dt$, which obviously is an
isomorphism between the first and second summands. Therefore, the $E^2$ term
and $HH_*(\D)$ vanish.

Of course, this criterion rarely applies in practice. In fact, the only
nontrivial case I know of are the Lefschetz fibrations obtained by Morsifying
the isolated critical point of the holomorphic functions $f(x,y,z) = xy +
z^{m+1}$. Like any other example obtained from singularity theory, the total
space $E$ is (deformation equivalent to) a ball, hence the symplectic homology
vanishes. The fibre $M$ is (a compact piece of) the so-called ALE space of type
$(A_m)$, and the standard basis of vanishing cycles is a chain of spheres, each
intersecting its neighbour in a single point. $\A$ is a kind of $(A_m)$ quiver
algebra, and $\B \iso \A \oplus \A^\vee[-2]$ (higher order compositions are
known to be irrelevant \cite{seidel-thomas99}).
\end{example}

Note that, even assuming that $\B$, hence also $\A$, is explicitly known,
$HH_*(\D)$ may not be computable in the abstract sense of the word, since
$\bar{C}^*$ is infinite-dimensional (this kind of situation is by no means new;
for instance, given a finitely presented group, the problem of computing its
group cohomology is generally unsolvable). However, the successive
approximations $\bar{C}^*/F^{p+1}\bar{C}^*$ are finite-dimensional, hence
accessible to computation, even though in practice, the steeply increasing
complexity tends to limit one to low values of $p$.

\begin{example}
Consider the Morsification of $f(x) = x^{m+1}$, with $m>1$. This can also be
described as a generic $(m+1)$-sheeted branched cover $\pi: E \rightarrow D$
where the total space is again $D$, so that symplectic homology vanishes (while
the argument from Example \ref{th:extension} does not apply in this case, it is
worth while noticing that the relevant natural transformation $\SS \rightarrow
[1]$ is at least nilpotent; this is for degree reasons, since the Serre functor
in $mod(\A)$ is known to have the property that $\SS^{m+1} \iso [m-1]$). We
take the simplest case $m = 2$, fire up our trusty laptop, and compute the
Betti numbers of $\bar{C}^*/F^{p+1}\bar{C}^*$, taking $K = \Q$:
\begin{equation}
\begin{array}{r|rrrrrrrrrrrr}
 & \hspace{6em}\text{degrees (homological)}\hspace{-14em} \\
 & -11 & -10 & -9 & -8 & -7 & -6 & -5 & -4 & -3 & -2 & -1 & 0 \\
 \hline
 p = 0 &&&&&&&&&&&& 2 \\
 p = 1 &&&&&&&&&& 3 & 1 & \\
 p = 2 &&&&&&&& 3 &&& 1 & \\
 p = 3 &&&&&& 3 &&& 1 &&& \\
 p = 4 &&&& 3 &&&&& 1 &&&\\
 p = 5 &&   3  &&&&&&& 1 &&&\\
 &&&& \cdots \hspace{-2em} \\
\end{array}
\end{equation}
This certainly seems compatible with the idea that the limit as $p \rightarrow
\infty$ vanishes.
\end{example}

\begin{example}
Take the fibration with fibre $M = DT^*\!S^{d-1}$ and two identical vanishing
cycles $L_0 = L_1 = S^{d-1}$. In this case, $E$ is deformation equivalent to
$DT^*\!S^d$, which allows one to use \eqref{eq:loop} to determine its
symplectic homology (there are a number of fibrations with similar properties,
obtained by complexifying real Morse functions \cite{johns06}; in this case,
the standard Morse function on $S^d$). We take $d = 2$, $K = \Q$, and compute:
\begin{equation}
\begin{array}{r|rrrrrrrrrr}
 & \hspace{6em}\text{degrees (homological)}\hspace{-14em} \\
 & -8 & -7 & -6 & -5 & -4 & -3 & -2 & -1 & 0 \\
 \hline
 p = 0 && &&&&&&& 2 \\
 p = 1 && &&&& 2 & 4 & 1 & 1  \\
 p = 2 && && 4 & 6 & 1 & 1 & 1 & 1 \\
 p = 3 && 6 & 8 & 1 & 1 & 1 & 1 & 1 & 1 \\
 &&&& \cdots \hspace{-2em} \\
 H^{-*}({\EuScript L}S^2;\Q) & \cdots & 1 & 1 & 1 & 1 & 1 & 1 & 1 & 1
\end{array}
\end{equation}
The coincidence becomes even more striking after passing to $\K = \Z/2$, where
we pick up extra torsion in cohomology (for the topological side, see
\cite{kuribayashi91}):
\begin{equation}
\begin{array}{r|rrrrrrrrrr}
 & \hspace{6em}\text{degrees (homological)}\hspace{-14em} \\
 & -8 & -7 & -6 & -5 & -4 & -3 & -2 & -1 & 0 \\
 \hline
 p = 0 && &&&&&&& 2 \\
 p = 1 && &&&& 2 & 4 & 1 & 1  \\
 p = 2 && && 4 & 6 & 2 & 2 & 1 & 1 \\
 p = 3 && 6 & 8 & 2 & 2 & 2 & 2 & 1 & 1 \\
 &&&& \cdots \hspace{-2em} \\
 H^{-*}({\EuScript L}S^2;\Z/2) & \cdots & 2 & 2 & 2 & 2 & 2 & 2 & 1 & 1
\end{array}
\end{equation}
\end{example}


\section{Periodic points\label{sec:periodic}}
We now try to probe the geometric meaning of $HH_*(\D)$, something
which is obviously central to a proper understanding of Conjecture
\ref{th:main}. The following discussion will unfortunately remain
somewhat incomplete, partly because it is based on other conjectural
material, but most importantly because it remains at the level of
$\CC$, which means of the starting term \eqref{eq:e1}. Moreover, we
will approach things in a somewhat roundabout way, involving a
detour through cyclic homology $HC_*(\CC)$. Recall that this sits in
a long exact sequence
\begin{equation} \label{eq:connes}
\cdots \rightarrow HH_*(\CC) \longrightarrow HC_*(\CC) \longrightarrow
HC_{*-2}(\CC) \rightarrow \cdots
\end{equation}
Moreover, since $\CC$ is unital and augmented,
\begin{equation} \label{eq:split-cyclic}
HC_*(\CC) = HC_*(R) \oplus \overline{HC}_*(\CC),
\end{equation}
where the second summand is reduced cyclic homology \cite[Section
2.2.13]{loday}. Note that all of this formalism is compatible with the direct
product decomposition by $t$-weights (where $HC_*(R)$ in
\eqref{eq:split-cyclic} is given weight $0$).

For simplicity, assume that $char(\K) = 0$. Then $\overline{HC}_*(\CC)$ can be
computed through the reduced Connes complex $\bar{C}^*_\lambda$, which is the
$t$-adic completion of
\begin{equation}
\bar{T}(\CC_+[1])^{cycl}[-1] = \bigoplus_{n=0}^\infty \big((\CC_+[1])^{\otimes
n+1,diag}\big)^{\Z/n+1}[-1].
\end{equation}
Here, the superscript $\Z/n+1$ denotes the coinvariant part for the action
which cyclically permutes factors (with signs). The differential is essentially
the same as for $HH_*(\CC)$. In the specific case of $\CC$, given any $x_n
\otimes \cdots \otimes x_0 \in (\CC_+[1])^{\otimes n+1, diag}$, one can permute
cyclically until the first element $x_n$ lies in $t^i\B$ for some $i>0$. This
means that the reduced Connes complex can be written as
\begin{equation} \label{eq:cyclic-pieces}
\prod_{k \geq 1} \Big(\!\!\!\prod_{i_1,\dots,i_k \geq 1}\!\! \big(t^{i_k}\B[1]
\otimes T(\A_+[1]) \otimes \cdots \otimes t^{i_1}\B[1] \otimes
T(\A_+[1])\big)^{diag}\Big)^{\Z/k}[-1].
\end{equation}
This time, the generator of $\Z/k$ acts by a permutation which moves the first
piece $t^{i_k}\B[1] \otimes T(\A_+[1])$ to the right end of the tensor product.
Note that the piece with $t$-weight one is simply $(t\B[1] \otimes
T(\A_+[1]))[-1]$, hence agrees with \eqref{eq:donaldson2} up to a shift by
$d+2$. Assuming Conjecture \ref{th:donaldson}, we can therefore interpret its
cohomology as $HF_{*+d+2}(\mu,id)$.

To extend this idea to higher powers of $t$, we move the base point $\ast$
slightly into the interior of $D$, and then take the $p$-fold cover $E^p$ of
$E$ branched along $M = E_\ast$. This is by itself the total space of a
Lefschetz fibration $\pi^p: E^p \rightarrow D^p \iso D$ with fibre $M$.
Concretely, given a basis of vanishing cycles for $\pi$, a $p$-fold repetition
of that list of Lagrangian spheres constitutes a basis for $\pi^p$. In
particular, the boundary monodromy is $\mu^p$. Now consider Conjecture
\ref{th:donaldson} for this new fibration: on one side we have
$HF_*(\mu^p,id)$, and on the other side (relying on the relation between
distinguished bases explained above) we have the cohomology of a complex of the
form
\begin{equation} \label{eq:hh-summand}
\bigoplus_{\substack{k \geq 1 \\ i_1+\cdots+i_k = p}} \K^{i_k} \otimes \big(
t^{i_k} \B[1] \otimes T(\A_+[1]) \otimes \cdots \otimes t^{i_1} \B[1] \otimes
T(\A_+[1])\big)^{diag}[d-1+2p].
\end{equation}
The branched cover case has one additional feature, which we need to take into
consideration. Conjugation by $\mu$ induces a $\Z/p$-action on $HF_*(\mu^p)$ as
well as its relative version $HF_*(\mu^p,id)$. In parallel, one can equip
\eqref{eq:hh-summand} with a $\Z/p$-action, as follows. Note that each summand
$(t^{i_k}\B[1] \otimes \cdots \otimes T(\A_+[1]))^{diag}$ occurs $i_k$ times,
which we distinguish by a label $1 \leq l \leq i_k$. The generator of the
action operates in the following way: if $l < i_k$, increase it by one and do
nothing; otherwise, apply a cyclic permutation as in \eqref{eq:cyclic-pieces},
and reset the label to zero. It seems plausible to expect the conjectural
isomorphism to be compatible with these two cyclic actions, and we will assume
without further ado that this is indeed the case. Then, the $\Z/p$-coinvariant
part of \eqref{eq:hh-summand} coincides with the $t$-weight $p$ piece of
\eqref{eq:cyclic-pieces} up to a shift by $d+2p$, generalizing the observation
previously made for $p = 1$. For each $p>0$, the weight $p$ piece of cyclic
homology would therefore be
\begin{equation}
HC_\ast(\CC,\CC)^p \iso \overline{HC}_\ast(\CC,\CC)^p \iso
HF_{\ast+d+2p}(\mu^p,id)^{\Z/p}.
\end{equation}
The outcome of these considerations is the following proposed geometric
interpretation:

\begin{conjecture} \label{th:periodic}
For each $p<0$, the column $E^1_{p\ast}$ of \eqref{eq:e1} sits in a long exact
sequence
\begin{equation} \label{eq:e1-explained}
\cdots \rightarrow E^1_{p\ast} \longrightarrow
HF_{\ast+d-p}(\mu^{-p},id)^{\Z/p} \longrightarrow
HF_{\ast+d-p-2}(\mu^{-p},id)^{\Z/p} \rightarrow \cdots
\end{equation}
\end{conjecture}

\begin{example}
Take a double branched cover $\pi: E \rightarrow D$ with an odd number $m>1$ of
branch points. In that case, the fibre $M$ consists of two points, and the
global monodromy $\mu$ exchanges those two points, or more precisely (taking
the grading into account) combines the exchange with a shift $[m]$. Moreover,
since there are no topological sections, the map \eqref{eq:sections} and its
analogues for $p>1$ must all vanish. Hence,
\begin{equation}
HF_*(\mu^p,id) = H_{*-1}(M;\K) \oplus HF_*(\mu^p) = \begin{cases} \K^2 & \ast =
1, \\ \K^2 & \ast = mp\text{, provided $p$ is even}, \\ 0 & \text{otherwise.}
\end{cases}
\end{equation}
The $\Z/p$-action is trivial on the first kind of generators but nontrivial on
the second one. Hence
\begin{equation}
E^1_{pq} = \begin{cases} \K^m & p = 0, q = 0, \\
\K^2 & \text{$p<0$, $q = p$ or $p+1$}, \\
\K & \text{$p<0$ even, $q = (1-m)p-1$ or $(1-m)p$.}
\end{cases}
\end{equation}
The only possibly nonzero differentials are $E^1_{00} = \K^m \rightarrow
E^1_{-10} = \K^2$ and $E^1_{pp} = \K^2 \rightarrow E^1_{p-1,p} = \K^2$ for
$p<0$. One can run a comparison argument with the case $m = 1$; the outcome
suggests that the first of these differentials is of rank $1$, while the rest
is of rank $2$. This would leave
\begin{equation} \label{eq:surface}
HH_*(\D,\D) = \begin{cases} \K^{m-1} & \ast = 0, \\
 \K & \ast = -1, \\
 \K & \ast = 2(m-2)-1, 4(m-2)-1, \dots, \\
 \K & \ast = 2(m-2), 4(m-2), \dots
\end{cases}
\end{equation}
On the other hand, the Reeb flow on $\partial E$ being just rotation, we have
$SH_*(E) = H_{*+1}(E;\K) \oplus H_{*+1+2(2-m)}(\partial E;\K) \oplus
H_{*+1+4(2-m)}(\partial E;\K) \oplus \cdots$, which agrees well with
\eqref{eq:surface}.
\end{example}

Let's draw some informal consequences of the discussion above. According to a
suitable version the standard definition, the chain complex underlying
$SH_*(E)$, for any Liouville domain $E^{2d}$, has one generator (of degree
$i-d$) for each critical point (of Morse index $i$) of an exhausting Morse
function on $int(E)$, and another pair of generators (in adjacent degrees,
governed by the Conley-Zehnder index) for each periodic orbit of the Reeb flow
on $\partial E$ (multiples of orbits are counted separately). In our case,
consider the boundary of $E$ (after rounding off corners) as a manifold having
an open book decomposition with page $M$ and monodromy $\mu$. In this
description, there is an obvious vector field transverse to the pages, and
periodic orbits of that vector field (excluding the constant ones inside the
spine) correspond to fixed points of $\mu^p$, for any $p$, mod the
$\Z/p$-action by $\mu$. One can associate to this open book decomposition a
contact structure (which is indeed the one we have been considering on
$\partial E$ all along), whose Reeb vector field is a small perturbation of the
previously mentioned transverse vector field. It follows that in a rough and
ready count, Conjecture \ref{th:periodic} accounts for most of the generators
needed to build $SH_*(E)$, the exceptions being the following: the generators
coming from the Morse function on $int(M)$; infinitely many copies of $H_*(M)$
(which make up the difference between $HF_*(\mu^p)$ and its relative version);
and finally, those periodic Reeb orbits which are located near the spine. Note
that these missing pieces are precisely those that enter into the definition of
$SH_*(D \times M)$, which is known to be zero. Therefore, it presumably makes
sense to think of $HH_*(\D)$ as an algebraic model for a relative symplectic
homology group $SH_*(E,D \times M)$, which then of course coincides with
$SH_*(E)$.

\section{Iterating the Serre functor\label{sec:iterate}}

According to the general theory, an $A_\infty$-module over $\D$ is a graded
$R$-module $\MM$ together with a structure map
\begin{equation} \label{eq:module}
\mu_\MM^*: \MM \otimes T(\D_+[1]) \longrightarrow \MM[1],
\end{equation}
whose components are the differential $\mu_\MM^1$ and the higher order (module
structure) terms. However, the proper definition also needs to take the
$t$-adic topology into account, which means that $\MM$ should carry a complete
decreasing filtration such that \eqref{eq:module} becomes a continuous map,
hence extends to the completed tensor product. For our present purpose, it will
be sufficient to consider modules for which this filtration is trivial, in
which case the condition says that \eqref{eq:module} should vanish on the part
of $\MM \otimes T(\D_+[1])$ where the second factor has sufficiently high
$t$-weight. Given two such modules $\MM_0$ and $\MM_1$, an element of
$hom(\MM_0,\MM_1)$ of degree $k$ is a map
\begin{equation}
\phi: \MM_0 \otimes T(\D_+[1]) \longrightarrow \MM_1[k],
\end{equation}
subject to the same vanishing condition as before. Finally, we impose an
additional condition that all $\MM$ should be finite as $R$-modules, and denote
by $D(\D)$ the resulting differential graded category. The truncation map $q:
\D \rightarrow \A$, defined by setting $t = 0$, gives rise to a pullback
functor $q^*: D(\A) \rightarrow D(\D)$. Concerning the behaviour of this, we
make the following

\begin{conjecture} \label{th:iterate}
For all $\MM_0,\MM_1 \in D(\A)$, there is an isomorphism
\begin{equation} \label{eq:mod-pullback}
H(hom_{D(\D)}(q^*\MM_0,q^*\MM_1)) \iso \underrightarrow{\lim}_p\;
H(hom_{D(\A)}(\SS^p(\MM_0),\MM_1)[dp]),
\end{equation}
where $\SS$ is the Serre functor of $\A$, and the connecting map in the
directed system on the right hand side is the one induced by \eqref{eq:alg-tr}.
\end{conjecture}

Informally, one can interpret this as saying that $D(\D)$ is produced from
$D(\A)$ through a process which tries to turn the natural transformation $\SS
\rightarrow id$ into an isomorphism. Note that, due to the uniqueness of Serre
functors, each group in the direct system on the right hand side is defined
entirely in terms of $\A$ (however, the maps connecting the groups do depend on
$\B$). We should also point out that Conjecture \ref{th:iterate} is purely
algebraic, and can be studied independently of its geometric motivation.

Let's try to see how the isomorphism \eqref{eq:mod-pullback} might possibly
come about, in the simplest example where both $\MM_0 = \MM_1 = R$ are the
simple module. Take the reduced bar construction $\bar{B}^* = \bar{B}^*(\D)$:
this is the $t$-adic completion of $T(\D_+[1])$, with a differential which
consists of \eqref{eq:d1} plus insertion of $\mu^0$ in all possible places,
like \eqref{eq:d3} but allowing $i = n+1$ as well. A simple calculation shows
that
\begin{equation}
 \bar{B}^* = hom_{D(\D)}(q^*R,q^*R)^\vee
\end{equation}
is the dual of the complex on the left hand side of \eqref{eq:mod-pullback}. In
analogy with \eqref{eq:pieces}, $\bar{B}^*$ can be written as the direct
product of pieces
\begin{equation} \label{eq:b-pieces}
T(\A_+[1]) \otimes t^{i_k}\B[1] \otimes \cdots \otimes t^{i_1}\B[1] \otimes
T(\A_+[1]).
\end{equation}
It comes with a complete $t$-adic filtration $F^*$, and we will be interested
in the finite approximations $\bar{B}^*/F^{p+1}$. For $p = 0$, all that remains
is $T(\A_+[1])$, so by analogy with the previous computation,
\begin{equation}
 \bar{B}^*/F^1 = hom_{D(\A)}(R,R)^\vee.
\end{equation}
In the next case $p = 1$, $\bar{B}^*/F^2$ contains a subcomplex of the form
\begin{equation} \label{eq:insert-a}
 Cone\big(T(\A_+[1])) \xrightarrow{m} T(\A_+[1]) \otimes t\A[1] \otimes
 T(\A_+[1]),
\end{equation}
where the map $m$ consists of insering $te_k$ in all possible positions.
Consider for a moment the classical situation, when $A$ is just an algebra over
$R$. In that case, $A \otimes T(A_+[1])$ is the standard bar resolution of the
$A$-module $R$, and $T(A_+[1]) \otimes tA[1] \otimes T(A_+[1])$ would therefore
be quasi-isomorphic to $R[-1] \otimes T(A_+[1]) = T(A_+[1])[-1]$, with the
quasi-isomorphism given precisely by $m$. Using a suitable filtration, this
argument carries over to the $A_\infty$-case, which shows that
\eqref{eq:insert-a} is acyclic. The quotient of $\bar{B}^*/F^2$ by
\eqref{eq:insert-a} is
\begin{equation} \label{eq:one-q}
T(\A_+[1]) \otimes t\QQ[1] \otimes T(\A_+[1]),
\end{equation}
where $\QQ$ is as in \eqref{eq:short}. Recall from \eqref{eq:derived-tensor}
that $T(\A_+[1]) \otimes \QQ = R \otimes T(\A_+[1]) \otimes \QQ$ is the tensor
product $R \otimes_\A \QQ$. Hence, \eqref{eq:one-q} is the dual of $hom(R
\otimes_\A \QQ[-1],R)$, and since $\QQ \iso \A^\vee[1-d]$, we conclude that
there is a quasi-isomorphism
\begin{equation}
 \bar{B}^*/F^2 \htp \big(hom_{D(\A)}(\SS(R),R)[d]\big)^\vee.
\end{equation}
After dualizing again, one sees that these are indeed the first two terms in
the direct system \eqref{eq:mod-pullback}. The behaviour for general $p$ is not
as easy to understand, but the following considerations might be helpful.
Suppose first that $\B = \A \oplus \QQ$ is the extension algebra constructed
from $\A$ and $\QQ$. In that case, $\bar{B}^*$ carries a natural decomposition
given by splitting each $\B$ in \eqref{eq:b-pieces}, and counting the number of
$\QQ$ factors. The part with the most such factors is
\begin{equation}
T(\A_+[1]) \otimes t\QQ \otimes T(\A_+[1]) \cdots \otimes t\QQ \otimes
T(\A_+[1]),
\end{equation}
with $p$ factors of $t\QQ$. In analogy with our previous analysis of
\eqref{eq:insert-a}, it is vaguely plausible to think that all the other pieces
might be acyclic. By using a suitable filtration, the result would then carry
over to general $\B$, yielding a quasi-isomorphism
\begin{equation}
\bar{B}^*/F^{p+1} \htp \big(hom_{D(\A)}(\SS^p(R),R)[dp]\big)^\vee.
\end{equation}
This, as should be clear from our formulation, is just speculation; rather than
pursuing it further, we finish our discussion by taking a brief look at the
geometric meaning of Conjecture \ref{th:iterate} in the original context of
Lefschetz fibrations.

Let's consider only objects of $D(\A)$ which correspond to actual Lefschetz
thimbles in $\F(\pi)$. Then, assuming Conjecture \ref{th:agree}, the right hand
side of \eqref{eq:mod-pullback} can be written as
\begin{equation} \label{eq:sigma-twist}
\underrightarrow{\lim}_p \, HF^*(\sigma^p(\Delta_0),\Delta_1).
\end{equation}
After smoothing the corners of $E$, this is apparently the same as taking
$\Delta_0$ and moving its boundary by the Reeb flow for increasingly large
times, then extending the isotopy to the interior and considering the resulting
sequence of Floer cohomology groups, in analogy with the definition of
symplectic homology itself. With this in mind, it seems even possible that
$D(\D)$ is itself an invariant of $E$.

\begin{example}
Suppose that $E$, after rounding off corners, is of the form $DT^*N$, and that
we have a Lefschetz thimble $\Delta$ which turns out to be isotopic to a
cotangent fibre. Set $\Delta_0 = \Delta_1 = \Delta$ in \eqref{eq:sigma-twist};
in analogy with \eqref{eq:loop}, it seems natural to expect the direct limit to
be the homology of the ordinary (based) loop space, $H_*(\Omega N;\K)$. This is
relevant to the mirror symmetry situation from Example \ref{th:mirror}, since
there $E$ is deformation equivalent to $DT^*(T^d)$, which would yield
\begin{equation} \label{eq:laurent}
H_*(\Omega T^d;\C) = \C[t_1^{\pm 1},\dots,t_d^{\pm}].
\end{equation}
To see how this fits in with the mirror symmetry prediction, let $U = X
\setminus Y \iso (\C^*)^d$ be the open torus orbit in $X$. In $D^b(X)$,
consider the natural transformation $S[-d] \rightarrow id$ given by $s$. For
any two objects $F_0,F_1$,
\begin{equation} \label{eq:complement}
\underrightarrow{\lim}_p \, Hom_{D^b(X)}(S^p(F_0),F_1)[dp] \iso
Hom_{D^b(U)}(F_0|U,F_1|U),
\end{equation}
since the direct limit just amounts to allowing poles of increasingly high
order along $Y$ (the most convenient way is to prove \eqref{eq:complement}
first for vector bundles, and then to argue through exact triangles). The
structure sheaf ${\mathcal O}_X$ is expected to correspond to a Lefschetz
thimble $\Delta$ which essentially is a cotangent fibre, see for instance
\cite{abouzaid05}; and indeed, the right hand of \eqref{eq:complement} for $F_0
= F_1 = {\mathcal O}_X$ is the affine coordinate ring $\C[U]$, hence isomorphic
to \eqref{eq:laurent}.
\end{example}



\bibliographystyle{plain}

\begin{thebibliography}{10}

\bibitem{abouzaid05}
M.~Abouzaid.
\newblock Homogeneous coordinate rings and mirror symmetry for toric varieties.
\newblock {\em Geometry and Topology}, 10:1097--1156, 2006.

\bibitem{bondal89}
A.~I. Bondal.
\newblock Representations of associative algebras and coherent sheaves.
\newblock {\em Math. USSR Izvestiya}, 34:23--42, 1990.

\bibitem{cieliebak02}
K.~Cieliebak.
\newblock Handle attaching in symplectic homology and the chord conjecture.
\newblock {\em J. Eur. Math. Soc.}, 4(2):115--142, 2002.

\bibitem{cieliebak-floer-hofer95}
K.~Cieliebak, A.~Floer, and H.~Hofer.
\newblock Symplectic homology {II}: a general construction.
\newblock {\em Math. Z.}, 218:103--122, 1995.

\bibitem{costello05}
K.~Costello.
\newblock Topological conformal field theories and {C}alabi-{Y}au categories.
\newblock Preprint math.QA/0412149.

\bibitem{donaldson02}
S.~K. Donaldson.
\newblock Lefschetz pencils on symplectic manifolds.
\newblock {\em J. Differential Geom.}, 53:205--236, 1999.

\bibitem{gorodentsev-kuleshov04}
A.~Gorodentsev and S.~Kuleshov.
\newblock Helix theory.
\newblock {\em Mosc. Math. J.}, 4:377--440, 2004.

\bibitem{johns06}
J.~Johns.
\newblock PhD thesis, University of Chicago, 2006.

\bibitem{keller}
B.~Keller.
\newblock Introduction to {$A$}-infinity algebras and modules.
\newblock {\em Homology Homotopy Appl.}, 3:1--35, 2001.

\bibitem{khovanov-seidel98}
M.~Khovanov and P.~Seidel.
\newblock Quivers, {F}loer cohomology, and braid group actions.
\newblock {\em J. Amer. Math. Soc.}, 15:203--271, 2002.

\bibitem{kuribayashi91}
K.~Kuribayashi.
\newblock On the mod {$p$} cohomology of the spaces of free loops on the
  {G}rassmann and {S}tiefel manifolds.
\newblock {\em J. Math. Soc. Japan}, 43:331--346, 1991.

\bibitem{lefevre}
K.~Lef{\`e}vre-Hasegawa. 
\newblock Sur les $A_\infty$-cat{\'e}gories.
\newblock Thesis, Paris 7 University, 2002.

\bibitem{loday}
J.-L.~Loday.
\newblock {\em Cyclic homology}.
\newblock Springer, 2nd edition, 1997.

\bibitem{oancea04b}
A.~Oancea.
\newblock A survey of {F}loer homology for manifolds with contact type boundary
  or symplectic homology.
\newblock {\em Ensaios Mat.}, 7:51--91, 2004.

\bibitem{oancea04}
A.~Oancea.
\newblock The {K}unneth formula in {F}loer homology for manifolds with
  restricted contact type boundary.
\newblock {\em Math. Ann.}, 334:65--89, 2006.

\bibitem{seidel00b}
P.~Seidel.
\newblock More about vanishing cycles and mutation.
\newblock In K.~Fukaya, Y.-G. Oh, K.~Ono, and G.~Tian, editors, {\em
  {S}ymplectic {G}eometry and {M}irror {S}ymmetry ({P}roceedings of the 4th
  {KIAS} Annual International Conference)}, pages 429--465. World Scientific,
  2001.

\bibitem{exact}
P.~Seidel.
\newblock A long exact sequence for symplectic {F}loer cohomology.
\newblock {\em Topology}, 42:1003-1063, 2003.

\bibitem{seidel04}
P.~Seidel.
\newblock {F}ukaya categories and {P}icard-{L}efschetz theory.
\newblock European Math. Soc., in press.

\bibitem{seidel-thomas99}
P.~Seidel and R.~Thomas.
\newblock Braid group actions on derived categories of coherent sheaves.
\newblock {\em Duke Math. J.}, 108:37--108, 2001.

\bibitem{tradler01}
T.~Tradler.
\newblock Infinity-inner-products on {A}-infinity-algebras.
\newblock Preprint math.\-AT/0108027.

\bibitem{viterbo97a}
C.~Viterbo.
\newblock Functors and computations in {F}loer homology with applications,
  {P}art {I}.
\newblock {\em Geom. Funct. Anal.}, 9:985--1033, 1999.

\bibitem{weber06}
J.~Weber.
\newblock Three approaches towards {F}loer homology of cotangent bundles.
\newblock {\em J. Symplectic Geom.}, 3(4):671--701, 2005.
\end{thebibliography}

\end{document}